\documentclass[10pt]{amsart}
\usepackage{amsmath,amssymb}

 \newtheorem{prop}{Proposition}
 \newtheorem{lem}{Lemma}
 \newtheorem{defi}{Definition}
 \newtheorem{cor}{Corollary}
  \newtheorem{theo}{Theorem}

\newcommand\s{\vee}
\newcommand\w{\wedge}

 \newcommand\p{^\perp}

\newcommand\pr{^{\prime}}

\newcommand\ep{^{-1}}
\newcommand\di{\displaystyle}
\begin{document}

 \title{Automorphisms of an orthomodular poset of projections}
\author{Georges Chevalier}
  \address{Georges Chevalier,
 Institut Girard Desargues,
 UMR 5028, Math\'{e}matiques,
 Universit\'{e} Claude-Bernard LYON 1,
 43, bd. du 11 novembre 1918,
 69622 Villeurbanne cedex,
 France.
  e-mail : chevalier@igd.univ-lyon1.fr}

\begin{abstract} By using a lattice characterization of continuous
projections defined on a topological vector space $E$ arising from
a dual pair, we determine the  automorphism group  of their
orthomodular poset $Proj(E)$ by means of automorphisms and
anti-automorphisms of the lattice $L$ of all closed subspaces of
$E$. A connection between the    automorphism group of the ring of
all continuous linear mappings defined on $E$ and the automorphism
group  of the orthoposet  $Proj(E)$ is established.
\end{abstract}

 \maketitle
\section{Introduction}
In a vector space $E$, there exists a natural correspondence
between projections and pairs of subspaces : to every projection
$p$ is associated the pair $(Im\, p\, ,\,Ker\, p)$ of subspaces.
If $E$ is a topological vector space and $p$ a continuous
projection then $(Im\, p\, ,\,Ker\, p)$ is a pair of closed
subspaces and $ Im\,p+Ker\,p $ is a topological direct sum. In a
previous paper (\cite{CHE}), we introduced the projection poset
$P(L)$ of a lattice $L$ satisfying some  properties of lattices of
closed subspaces. We proved that if $L$  is   the lattice   of all
closed subspaces of a topological vector space $E$ arising from a
dual pair, then $P(L)$ is isomorphic to the poset of continuous
projections defined on $E$ (Theorem 1 and 2 of \cite{CHE}).  By
using this isomorphism, we determined   the automorphism group of
a poset of continuous projections   by means of automorphisms and
anti-automorphisms of the lattice $L$ (Theorem 3 of \cite{CHE}).

This paper continues \cite{CHE}. In the first part,  Theorem 3 is
improved, a restrictive hypothesis is removed   and its setting is
extended to some incomplete lattices.

In the second part of the paper, we prove a continuous form of the
first fundamental theorem of projective geometry. This result
allows us to relate the     automorphism group of the orthomodular
poset of continuous projections defined on a topological vector
space $E$  with the automorphism group of the ring of continuous
linear mappings defined on  $E$.

Information about the lattice concepts used in this paper may be
found in \cite{MAE},  and \cite{KO} or \cite{SCH} are good
references for topological vector spaces.

 \section{  The orthomodular poset of projections of a symmetric
lattice.}

 In this section, we recall some definitions and results from
 \cite{CHE} where the reader is referred to for more information.

  In a lattice $L$, $(a,b)\in L^2$ is a   modular pair, written
 $(a,b)M$, if $(x\vee a)\wedge b= x\vee(a\wedge b)$ for every
$x\leq b$.
 The pair
$(a,b)$ is a   dual modular pair, written $(a,b)M^*$, if $(a,b)M$
holds in the dual lattice $L^*$ of $L$ and the lattice $L$ is said
to be a   symmetric lattice  if $(a,b)M$ implies $(b,a)M$ and
$(a,b)M^*$ implies $(b,a)M^*$.

Our purpose in the following definition of the projection poset of
a lattice $L$ is to obtain, when $L$ is the lattice of all closed
subspaces of a topological vector space $E$, a poset defined in a
algebraic setting and isomorphic to the poset of all continuous
linear
  projections defined on $E$. See \cite{CHE} for a
discussion about the motivation of this definition.

\begin{defi} Let $L$ be a symmetric lattice. The projection poset $P(L)$
 of $L$ is the following subset
 of the direct product $L\times L$ :
 $$  P(L)=\{(a,b)\in L\times L\;|\; a\vee b = 1,\; a\wedge b =0,\;
 (a,b)M,\;(a,b)M^*\}  $$
  \end{defi}
For a projection $p=(a,b)\in P(L)$, $a$ is called the image of $p$
and $b$ its kernel.

  If $(a,b)$ is a projection  of a symmetric lattice $L$ then
$(b,a)$ is also a projection and we write
   $(b,a) =(a,b)\p$.

  \begin{prop}(\cite{MU},  \cite{HA},  \cite{CHE})
 Let $L$ be is a symmetric lattice with 0 and 1. If $P(L)$ is ordered by the
 restriction $\leq$ of the order relation on $L\times L^*$
 then
 $(P(L),\leq,\perp)$
  is an orthomodular poset (abbreviated OMP). If $L$ possesses a
 structure of OMP then this OMP is naturally
  isomorphic to a
 suborthomodular poset of $P(L)$.
 \end{prop}

   An AC-lattice is an atomistic lattice with the covering property :   if $p$
 is an atom and $a\w p=0$ then $a\lessdot a\s p$ that is $a \leq x \leq a\s p$
  implies $a=x $ or $a\s p =x$.
   \medskip

 If $L$ and $L^*$ are AC-lattices, $L$ is called a   DAC-lattice.
 Any DAC-lattice is  symmetric and finite-modular (\cite{MAE}, Theorem 27.5).
Irreducible complete DAC-lattices of length $\geq 4$ are
 representable by lattices of closed subspaces   and many lattices
 of subspaces are DAC-lattices. We will now specify this last
 assertion.

 Let $K$   be a field, $E$ a left vector spaces over $K$, $F$ a
 right vector space over $K$.
 If there exists   a nondegenerate   bilinear form $
 \mathcal{B}$ on $E\times
 F$, we say that $(E, F)$ is a pair of dual spaces. For example, if $E$ is a locally
  convex space and $E\pr$ its topological dual space then $(E,E\pr)$ is naturally a pair
   of dual spaces with $\mathcal{B}(x,y)= y(x)$ (\cite{KO}, page 234).

     For a
 subspace $A$ of $E$, we put $$A\p = \{y\in F\;|\;\mathcal{B}(x,y) =0  \mathrm{\;
 for\;
 every}\; x \in A\}.$$
 Similarly, let

$$B\p = \{x\in E\;|\; \mathcal{B}  (x,y) =0  \mathrm{\;
 for\;
 every}\; y \in B\}$$
 for every subspace $B$ of $F$.
 A subspace $A$ of $E$ is called   $F$-closed  if $A= A^{\perp \perp}$
 and the  set of all $F$-closed subspaces, denoted by $L_F(E)$ and ordered by set-inclusion, is a complete
 irreducible DAC-lattice. Conversely,
 for any irreducible complete DAC-lattice $L$ of length $\geq 4$,
 there exists a pair $(E,F)$ of dual spaces such that $L$ is
 isomorphic to the lattice of all $F$-closed subspaces of $E$ (\cite{MAE}, Theorem 33.7).

 The set of all $E$-closed subspaces of $F$ is similarly defined and
 is also
 a DAC-lattice.

Let $(E,F)$ be a pair of dual spaces. The     linear weak topology
on $E$, denoted by $\sigma(E,F)$, is the linear topology defined
by taking $\{G\p\;|\; G\subset F,\;\dim\,G <\infty\}$ as a basis
of neighbourhoods of 0. If $F$ is interpreted as a subspace of the
algebraic dual of $E$ then a subbasis of neighbourhoods of 0
consists   of  kernels of elements of $F$.

 The  linear weak topology on $F$,
noted $\sigma(F,E)$, is defined in the same way. The space $F$ can
be interpreted as  the topological dual of $E$ for the
$\sigma(E,F)$ topology and $E$ as the topological dual of $F$ for
the $\sigma(F,E)$ topology. Equipped with their linear weak
topologies, $E$ and $F$ are topological vector spaces
(\cite{KO},{\S} 10.3) if  the topology on $K$ is   discrete.

 Moreover, for a subspace
$G\subset E$, we have $\overline{G}=G^{\perp\perp}$ and thus a
closed subspace in $E$ is an unambiguous notion.

The following theorem shows that our definition of a projection
poset of a lattice  is appropriate for our purpose.

 \begin{theo}\label{continuous}(\cite{CHE}) Let $L$ be a complete irreducible
  DAC-lattice. If $L$ is representable as
the lattice $ \mathcal{L}$ of all  $F$-closed subspaces of a pair
of dual spaces $(E,F)$   then the projection orthoposets $P(L)$
and $P(\mathcal{L})$ are isomorphic and  the correspondence $p
\mapsto (Im\,p,Ker\,p)$ is an isomorphism between the orthomodular
poset of $\sigma(E,F)$-continuous linear projections defined on
$E$ onto
  $P(\mathcal{L})$
\end{theo}

The   linear weak topology seems to be a poor topology. However, a
linear mapping $f$, defined on a locally convex space $E$, is
weakly continuous if and only $f$ is continuous for the linear
weak topology $\sigma(E,E')$ (\cite{KO},20.4) and so we obtain the
following
 consequences
  of Theorem \ref{continuous}.
\begin{cor}\label{local}(\cite{CHE}) Let $E$ be a locally convex space and $L$ its lattice
of all closed subspaces. The projection orthomodular poset $P(L)$
is isomorphic to the poset of weakly continuous linear projections
defined on  $E$.
\end{cor}
\begin{cor}(\cite{CHE}) If $H$is a  Hilbert space (more generally, a metrizable space)
and $L$ its lattice of closed subspaces then the projection
orthomodular poset $P(L)$ is isomorphic to the orthoposet of
  continuous linear projections defined on  $H$.
  \end{cor}

  \section{Automorphisms of an orthomodular poset of projections}
  The main result of  \cite{CHE}  is a generalization of a theorem of  \cite{OV}  and gives a description of
  automorphisms of a projection orthoposet $P(L)$ by means of
  automorphisms and anti-automorphisms of the lattice $L$ when $L$ is a
  complete DAC-lattice satisfying the condition

  \begin{center} for every $a\in L$ there exists $b\in L$ such
  that $(a,b)\in P(L)$   (C).
  \end{center}
  Moreover, it is proved in  \cite{CHE}  that there are   exactly two kinds of
  automorphisms on an orthoposet of projections : the so-called even automorphisms  which
  transform projections with the same image into projections with
  the same image and the odd automorphisms which transform
  projections with the same image into projections with the same
  kernel. This fact generalizes a theorem of \cite{OV}.

  In this section, we will improve on the main result of \cite{CHE} by
  removing the restriction condition (C)
    and by extending its setting to certain
  incomplete lattices.

  Let us say that an  irreducible DAC-lattice $L$ is a   G-lattice  if $L$ is complete
  or if $L$ is modular and complemented. Typical examples of G-lattices are obtained by considering
     a Hilbert space $H$: the lattice of
  all closed
  subspaces of $H$ is a G-lattice as a complete irreducible
  DAC-lattice  and its sublattice of   finite or cofinite dimensional
  elements is a G-lattice as a complemented modular irreducible
  DAC-lattice.
  Irreducible DAC-lattices of length $\geq  4$ which are either
  complete or modular and complemented share the following
  properties :
\begin{itemize}
  \item Every atom has more than one complement;
  \item  If $a\lessdot b$ then there exist different atoms $p_1$
  and $p_2$ such that $a\vee p_1= a\vee p_2 = b$;
  \item Two different atoms have a common complement.
\end{itemize}
By using these  facts,  all the results preceding Theorem 3 of
\cite{CHE}, proved for irreducible complete DAC-lattices,  extend
to G-lattices and an improved version of Theorem 3 is as follows.

  \medskip
  \begin{theo}\label{main} Let $L$ be a G-lattice of length $\geq 4$.
For every   automorphism $\phi$ of the poset $P(L)$ there exists
\begin{enumerate}
  \item  an automorphism $f$ of the lattice $L$ such that $\phi ((a,b))=
  (f(a),f(b))$, $(a,b) \in P(L)$, if $\phi $ is even,
  \item an anti-automorphism $g$ of the lattice $L$ such that $\phi ((a,b))=
  (g(b),g(a))$, $(a,b) \in P(L)$, if $\phi $ is odd.

\end{enumerate}
Conversely, if $f$ is an automorphism of $L$ then $\phi : P(L)
\mapsto L\times L^* $ defined by $\phi ((a,b))=
  (f(a),f(b))$ is an even automorphism of $P(L)$ and if $g$ is an
  anti-automorphism of $L$ then $\psi : P(L) \mapsto L\times L^*$ defined
  by $\psi ((a,b))=
  (g(b),g(a))$ is an odd automorphism of $P(L)$.
\end{theo}
\noindent Proof. Fist we recall some notations from \cite{CHE}. In
a DAC-lattice $L$, $At(L)$ denotes the set of all atoms, $At^*(L)$
is the set of all coatoms, and $ \mathcal{F}(L)$ is the G-lattice
of all finite or cofinite elements of $L$. By $P_1(L)$  we mean
the set of all atoms of the projection poset $P(L)$.

 Let us denote by $L^+$ the lattice $L$ if $L$ is
an irreducible complemented modular DAC-lattice and the lattice $
\mathcal{F}(L)$  if $L$ is a complete irreducible DAC-lattice. In
the two cases, $L^+$ is an irreducible complemented modular
DAC-lattice  and the restriction of $\phi$ to $P(L^+)$ is an
automorphism.

 Assume that $\phi$ is even. By Proposition
 8 of \cite{CHE}, there exist  two bijections $f_1:At(L^+)\mapsto
At(L^+)$ and $f_2:At^*(L^+)\mapsto At^*(L^+)$ such that, for every
$(p,q) \in P_1(L^+)$, $\phi ((p,q)) = (f_1(p), f_2(q))$. Let $a\in
L^+$, $a\neq 0$. There exists $b\in L^+$ such that $(a,b)\in
P(L^+)$. For any atom $p\leq a$ there exists a coatom $q$ with
$(p,q)\leq (a,b)$. Thus $\phi((p,q)) = (f_1(p),f_2(q))\leq
\phi(a,b)$. If $\phi((a,b))=(c,d)$ then $$f_1(\{p\in At(L^+)\;|\;
p\leq a\}) \subset \{p\in At(L^+)\;|\; p\leq c\}.$$ By using
$\phi\ep$, we have   $$f_1(\{p\in At(L^+)\;|\; p\leq a\}) = \{p\in
At(L^+)\;|\; p\leq c\} $$ and so  Proposition 9 of \cite{CHE}
implies that  $f_1$ can be extended to an automorphism $
\overline{f_1}$ of the lattice $L^+$. Similarly, $f_2$ has an
extension, $ \overline{f_2}$.

The correspondence $(a,b)\in P(L^+) \mapsto ( \overline{f_1}(a),
\overline{f_2}(b))$ is an automorphism of the poset $P(L^+)$ which
agrees with $\phi$ on $P_1(L^+)$. As $P(L^+)$ is atomistic (Lemma
6 of \cite{CHE}) , for every $(a,b) \in P(L^+)$, we have
$\phi((a,b))=( \overline{f_1}(a), \overline{f_2}(b))$. This
equality is also true for $(a,b)=(0,1)$.

By Proposition 6 of \cite{CHE}, $\phi$ is also an automorphism of
the orthoposet $P(L^+)$ and  thus $\phi((a,b)\p)=  \phi
((a,b))\p$, that is $( \overline{f_1}(b), \overline{f_2}(a))=(
\overline{f_2}(b), \overline{f_1}(a))$
 and so $ \overline{f_1}= \overline{f_2}$.

 The proof is similar if $\phi $ is odd and is complete if $L$ is a complemented modular
 DAC-lattice. If $L$ is an irreducible complete DAC-lattice, a
 lemma is necessary.

 \begin{lem}\label{complete} Let $L$ be an irreducible complete
 DAC-lattice of length $\geq 4$. Any
   automorphism of the lattice $ \mathcal{F}(L)$ extends to an
   automorphism of $L$.
 \end{lem}
\noindent Proof. Let $(E,F)$ be a pair of dual spaces such that
$L$ is
 isomorphic to the lattice $L_F(E)$   of all $F$-closed subspaces of
 $E$. The lemma will be proved if any automorphism  $\psi$ of $
 \mathcal{F}(L_F(E))$ extends to an automorphism of $L_F(E)$.

   Define, for every subspace
  $N$ of $E $,   $\varphi(N) = \bigcup    \{ \psi(M)\;|\; M\subset
  N,\; \dim M <\infty\}$. It is clear that $\varphi(N)$ is a
  subspace of $E $.
   Let $X$ be a subspace of $E $ and
$N =\bigcup    \{ \psi\ep (M)\;|\; M\subset
  X,\; \dim M <\infty\}$. The set $N$ is a subspace of $E $ and
   we have   $\varphi(N)=\bigcup    \{ \psi(\psi\ep(M))\;|\; M\subset
  X,\; \dim M <\infty\}= X$.

  Let $M$, $N$  be two subspaces of $E $. If $M\subset N$ then
  $\varphi(M)\subset \varphi(N)$ and, for the converse, let $L$ be

  a subspace of $M$ with $\dim L=1$. We have $\psi(L)\subset
  \varphi(M)\subset \varphi(N)$ and if $0\ne x \in\psi( L)$ then  there exists a subspace $K
  \subset N$, $\dim K <\infty$, such that $ x\in  \psi(K) $. By $\dim \psi(L) =1$, we  have
  $\psi(L)\subset \psi(K)$ and therefore
    $L\subset K \subset N$. Finally, $ M \subset N $ and $\varphi$ is an
    automorphism
  of the lattice  of all subspaces of $E $.
   This automorphism extends
  $\psi$ since, for $M\in L_{F }(E)$, $\psi (M)$ and $\varphi(M)$ have
  the same finite dimensional subspaces.

  Let $M\in L_F(E)$. As  $L_F(E)$ is a DAC-lattice there exists a
  family $(H_{\alpha})$ of $F$-closed hyperplanes such that
  $M=\bigwedge H_{\alpha} = \bigcap H_{\alpha}$ and
   $$\varphi (M) = \varphi( \bigcap H_{\alpha}) =  \bigcap \varphi(H_{\alpha})= \bigcap \psi(H_{\alpha})=
   \bigwedge
  \psi(H_{\alpha}).$$
    Therefore $\varphi(M)$ is $F$-closed and,
   as $\varphi^{-1}(M)$ is also $F$-closed, $\varphi$ is an
   automorphism of $L_F(E)$ extending $\psi$.

   We return to the proof of the theorem. If $\phi$ is   an even
   automorphism of $P(L)$, $L$ an irreducible complete
   DAC-lattice,
   then $\phi$ is also an automorphism of $P(\mathcal{F}(L))$ and
   so there exists an automorphism $f$ of $ \mathcal{F}(L)$  such
   that $\phi(a,b) = (f(a),f(b))$ for any $(a,b) \in
   \mathcal{F}(L)$. By using the lemma, $f$ extends to an
   automorphism of the lattice $L$ and, as $P(L)$ is an atomistic
   lattice, $\phi(a,b) = (f(a),f(b))$ for any $(a,b) \in P(L)$.

   The proof is similar if $\phi$ is odd.
   \medskip

 For the converse, Proposition 5 of \cite{CHE} shows that it suffices to prove that
 $(a,b) \in P(L)$ implies $(f(a), f(b))\in P(L)$ for any automorphism
 $f$ and $(g(b),g(a))\in P(L)$ for   any anti-automorphism $g$.
  But  these implications are an easy consequence of the
  equivalence :
   $$(a,b)M\Leftrightarrow  \forall x\in L,\; ((x\w b)\s a)\w b=
(x\w b)\s (a\w b)$$ and $$(a,b)M^*\Leftrightarrow \forall x\in
L,\; ((x\s b)\w a)\s b= (x\s b)\w (a\s b).$$

  \section{More about automorphisms}

  By Theorem \ref{main}, the automorphisms of the projection poset
  $P(L)$ of a complete irreducible DAC-lattice $L$ of length $\geq 4$
  are determined
  by the automorphisms and the anti-automorphisms of the lattice
  $L$. As every complete irreducible DAC-lattice of length $\geq 4$ is the lattice of all
  closed subspaces of a pair of dual spaces, in this section we will investigate
 the  automorphism group of the lattice
  of closed subspaces.

  \subsection{A continuous form of the first  fundamental theorem  of projective
  geometry.}

  If $E_1$ and $E_2$ are vector spaces of dimensions at least 3 over the fields $K_1$ and
  $K_2$ then, by the first  fundamental theorem of projective geometry (\cite{BA},
  page 44 or \cite{VA}, page 21), the lattices of all subspaces of
  $E_1$ and $E_2$ are isomorphic if and only if $K_1$ and $K_2$
  are isomorphic fields and $E_1$ and $E_2$ have the same
  dimension. Moreover, if $\psi$ is an isomorphism from the lattice
  of all subspaces of $E_1$ onto the lattice of all subspaces of
  $E_2$ then there exists a semi-linear bijection $s:E_1\mapsto
  E_2$ such that, for every subspace $M\subset E_1$,
  $\psi(M)=s(M)$. Conversely, every semi-linear bijection of $E_1$
  onto $E_2$ induces a lattice isomorphism.

  In the following proposition, we generalize
    a part of the previous result to lattices of closed subspaces.
  \begin{prop} Let $(E_1,F_1)$ and $(E_2,F_2)$ be two pairs  of dual spaces
 over the fields $K_1$ and $K_2$. If there exists an isomorphism
$\psi$ of the lattice $L_{F_1}(E_1)$ onto the lattice
$L_{F_2}(E_2)$ then $K_1$ and $K_2$ are isomorphic fields and
there exists a semi-linear bijection $s:E_1 \mapsto E_2$ such
that, for every $F_1$-closed subspace $M$ of $E_1$, $\psi(M) =
s(M)$
\end{prop}
Proof. The mapping $\psi$ is an order isomorphism of the poset of
all finite dimensional subspaces of $E_1$ onto the poset of all
finite dimensional subspaces of $E_2$.

  Define, for every subspace
  $N$ of $E_1$,   $\varphi(N) = \bigcup    \{ \psi(M)\;|\; M\subset
  N,\; \dim M <\infty\}$. By a proof similar to the proof of Lemma
  \ref{complete}, $\varphi$ is an isomorphism of the lattice  of all
  subspaces of $E_1$ onto the lattice of all subspaces of $E_2$
   which extends
  $\psi$. Thus, by the first  fundamental theorem    of projective
  geometry, the fields $K_1$ and $K_2$ are isomorphic and
  there exists a
semi-linear bijection $s:E_1 \mapsto E_2$ such that, for every
$F_1$-closed subspace $M$ of $E_1$, $\psi(M) = s(M)$

\medskip

\noindent{\bf Remark .}
 This proof is similar to the proof of
Lemma 1 of \cite{FIL} where the authors prove the same result for
complex normed spaces.

\medskip

In the case of lattices of all subspaces of vector spaces, any
semi-linear bijection induces a  lattice isomorphism. For lattices
of closed subspaces, only continuous semi-linear bijections are
allowed.

  \begin{prop}\label{hyper} Let $(E_1,F_1)$ and $(E_2,F_2)$ be two pairs
   of dual spaces over the same field. If  $E_1$ and $E_2$ are equipped, respectively, with the
   $\sigma(E_1,F_1)$-topology and the $\sigma(E_2,F_2)$-topology
   then,
    for every semi-linear bijection $s:E_1\mapsto E_2$, the following
  statements are equivalent.
\begin{enumerate}

  \item The bijection $s$ is  bicontinuous (i.e. both $s$ and $s^{-1}$ are continuous).
  \item $H\in L_{F_1}(E_1) \mapsto s(H)$ is a bijection from the set of all
  $F_1$-closed hyperplanes of $E_1$ onto the set of all $F_2$-closed hyperplanes of $E_2$.
  \item  $M\in L_{F_1}(E_1) \mapsto s(M)$ is an isomorphism from the
 lattice $L_{F_1}(E_1)$ onto $L_{F_2}(E_2)$.
\end{enumerate}

\end{prop}
Proof. $1)\Rightarrow 2)$. Since $s$ is a semi-linear bijection,
the correspondence $M \mapsto s(M)$ is an isomorphism of the
lattice of all subspaces of $E_1$ onto the lattice of all
subspaces of $E_2$ and maps bijectively   the sets of all
hyperplanes. If $H\subset E_2$ is an $F_2$-closed hyperplane then
$H$ is a neighbourhood of 0 for the $\sigma(E_2,F_2)$ topology.
Since $s$ is
 continuous,   there exists a finite dimensional subspace
  $G \subset F_1$ such that $G^\perp \subset s\ep (H)$.     As $s\ep (H)$ has a finite codimension
  in $G^\perp$, $s\ep (H)$ is closed (\cite{KO}, property (7), page
  87) and since $s\ep$
is also continuous, $H \mapsto s(H)$ is a bijection from the set
of all $F_1$-closed hyperplanes of $E_1$ onto the set of all
$F_2$-closed hyperplanes of $E_2$.\\
 $2)\Rightarrow 3)$. Let $M\in L_{F_1}(E_1)$. As $L_{F_1}(E_1)$ is a
DAC-lattice there exists a family $(H_\alpha)$ of $F_1$-closed
hyperplanes such that  $M = \bigwedge H_\alpha = \bigcap H_\alpha$
and therefore $s(M)  = s(\bigcap H_\alpha) = \bigcap s( H_\alpha)
=\bigwedge s( H_\alpha) $. Thus $s(M)\in L_{F_2}(E_2)$ and
$s(L_{F_1}(E_1)) \subset L_{F_2}(E_2)$. As $s\ep $ also satisfies
the statement $(2).$, $s(L_{F_1}(E_1)) = L_{F_2}(E_2)$ and $M\in
L_{F_1}(E_1) \mapsto s(M)$ is an isomorphism from the
  lattice $L_{F_1}(E_1)  $ onto $L_{F_2}(E_2)  $

  \noindent $3) \Rightarrow 1)$. This is clear  since the family of all
   closed hyperplanes is a 0-neighbourhood subbasis for the   linear weak topology.
\medskip

    \bigskip

    \begin{cor}
    Let $E_1$ and $E_2$ be   real metrizable locally convex spaces.
    \begin{enumerate} \item $E_1$ and $E_2$ are isomorphic if and
    only if their lattices of closed subspaces $ \mathcal{C}(E_1)$ and $ \mathcal{C}(E_2)$ are isomorphic.
    \item $\psi : \mathcal{C}(E_1)\mapsto \mathcal{C}(E_2)$ is a
    lattice isomorphism if and only if there exists a bicontinuous
    linear bijection $s:E_1 \mapsto E_2$ such that, for every $M\in
    \mathcal{C}(E_1)$, $\psi(M)=s(M)$.
    \end{enumerate}
    \end{cor}
Proof. As $E_1$ and $E_2$ are real vector spaces, semi-linear
bijections are simply linear bijections and, as $E_1$ and $E_2$
are metrizable locally convex   spaces, a linear mapping $s:E_1
\mapsto E_2$ is continuous if and only if $s$ a continuous mapping
for the linear weak topologies.
\medskip

\noindent{\bf Remark.}  This corollary is a generalization of the
following result of Mackey (\cite{MAC}): two real normed spaces
$X_1$ and $X_2$  are isomorphic if and only if there exists a
linear bijection $T : X_1\mapsto X_2$ which carries bijectively
closed hyperplanes of $X_1$ into closed hyperplanes of $X_2$; if
$T$ exists then $T$ is bicontinuous. This result is extended to
complex normed spaces in \cite{FIL} : if $\psi : \mathcal{C}(X)
\mapsto \mathcal{C}(Y)$ is an isomorphism of the lattices of
closed subspaces of infinite dimensional complex normed spaces $X$
and $Y$ then there exists a bicontinuous linear or conjugate
linear bijection $s: X\mapsto Y$ such that $\psi(M) = s(M)$ for
all $M\in  \mathcal{C}(X)$.  By using the   following theorem,
this last result is also  a consequence of Proposition
\ref{hyper}: if $s : X\mapsto Y$ is a bijective semi-linear
  transformation of infinite-dimensional complex normed spaces
  that carries closed hyperplanes to closed hyperplanes then $s$ is
  either linear or conjugate linear (\cite{KA} or \cite{FIL},
  Lemma 2).

\subsection{Automorphism group of a projection poset}
Let $(E,F)$ be a pair of dual spaces.

If $L$ is the DAC-lattice of all $F$-closed subspaces of $E$ then
Theorem \ref{main} and Proposition \ref{hyper} allow one to obtain
a description of       even automorphisms of the projection
lattice $P(L)$ by means of lattice automorphisms of $L$. They are
all the correspondences $(M,N)\in P(L) \mapsto (s(M),s(N))$ where
$s: E\mapsto E$ fulfills the equivalent conditions of Proposition
\ref{hyper}. In the language of rings, even automorphisms of the
OMP of all continuous linear projections are of the form $p\mapsto
sps\ep$ since, for a linear projection $p$, $s(Im\, p)=Im\,sps\ep$
and $s(Ker\,p) =Ker\,sps\ep$.

 If $g$ is an anti-automorphism of $L$ then any
anti-automorphism is of the form $fg$ where $f$ is some
automorphism. Thus the set of all anti-automorphisms is determined
by a particular anti-automorphism and the group of all
automorphisms but  it seems difficult to find conditions assuring
the existence of an anti-automorphism of $L=L_F(E)$. We will now
discuss this point.

A first case is well-known : if $L$ is the modular lattice of all
subspaces of an infinite-dimensional vector space $E$   then $L$
has no anti-automorphism (\cite{BA}, Self-duality theorem, page
97). This results is extended in \cite{OR}  to   infinite
dimensional projective geometries that are  irreducible,
complemented, modular, complete, atomic lattices of infinite
length. Such lattices can be represented as lattices of closed
subspaces of pairs of dual spaces of infinite dimension    in
which any sum of two closed subspaces is closed (\cite{OR}). Thus
the automorphism group of the OMP of projections defined on an
infinite-dimensional projective geometry $L$ is isomorphic to the
automorphism group of the lattice   $L$.

Now suppose that $E$ is finite dimensional. As $F$ is isomorphic
to a subspace of $E^*$ and $E$ to a subspace of $F^*$, we can
assume that $E=F$ with $dim \,E=n$. The existence of an
  anti-automorphism of L is equivalent to the existence
of an   anti-automorphism of the field $K$ :   if $\alpha$ is an
anti-automorphism of $K$ and if $(e_i)_{1\leq i\leq n }$ is a
basis of $E$ then   the $\alpha$-bilinear form $$(\di\sum x_ie_i
 ,\di\sum y_ie_i)  \mapsto \langle \di \sum x_ie_i,
\di \sum y_ie_i \rangle = \di\sum _{i=1}^n x_i\alpha(y_i)$$ is
non-degenerate and determines  an   anti-automorphism $g$ of $L$.

The anti-automorphism $\alpha$ is involutary if and only if $g$ is
involutary and, in this case, the group formed by the
automorphisms and the anti-automorphisms of $L$ is the semi-direct
product of the normal subgroup of all automorphisms and the
subgroup $\{1_E,g\}$. The automorphism group of $P(L)$ is the
semi-direct product of the normal subgroup of even automorphisms
and a two-element  subgroup $\{1_{P(L)}, \gamma\}$ where $\gamma$
is an involutary odd automorphism.

In the infinite dimensional case and if $E=F$ then a particular
anti-automorphism is given by $$X\in L_E(E) \mapsto X\p = \{x \in
E \;|\; \langle x,X\rangle =0\}$$ and the previous results allows
one to obtain all the automorphisms and all the anti-automorphisms
of $L_E(E)$ and thus to determine the automorphism group of its
projection lattice.

\medskip

\noindent{\bf Example.} If $H$ is a Hilbert space then the two
pairs of dual spaces $(H,H)$ and $(H,H\pr)$ coincide. The
correspondence $X\in L_H(H)\mapsto X^\perp$ is an involutary
anti-automorphism of the lattice of all closed subspaces of $H$
and $p\mapsto p^*$ is the corresponding  involutary odd
automorphism of the orthoposet of continuous linear projection
defined on $H$. By using the previous results, we find again the
main result of \cite{OV} : the automorphisms of the orthoposet
$proj (H)$ of all continuous linear projections of $H$ are of the
form $p\mapsto s\ep p s$ or $p\mapsto s\ep p^* s$ where $s$ is a
continuous linear bijection in the real case and a continuous
linear or conjugate linear bijection in the infinite dimensional
complex case. In the finite dimensional complex case, $s$ is only
a semi-linear bijection. The automorphism group of $Proj(H)$ is a
semi-direct product of the normal subgroup of even automorphisms
and a two-element subgroup.
\subsection{Application
to the ring of continuous linear mappings}

 The following proposition is proved in \cite{ JA}  (Isomorphism Theorem, page 79)  in the study of primitive rings
  having minimal right ideals. In \cite{FIL}, a different proof is given in the particular case of
infinite-dimensional complex normed linear  spaces. Here, we
generalize the latter proof    for two pairs of dual spaces and
obtain the  result   of \cite{JA}.
\begin{prop} Let $(E_1,F_1)$ and $(E_2,F_2)$ be two pairs  of dual spaces
(over the fields $K_1$ and $K_2$) and let us denote by $
\mathcal{B}(E_1)$ and $ \mathcal{B}(E_2)$ the rings of all
continuous linear mappings defined on $E_1$ and $E_2$ equipped
with their
 linear weak topologies. If there exists an isomorphism of rings,  $\Phi : \mathcal{B}(E_1) \mapsto
\mathcal{B}(E_2)$,    then $K_1$ and $K_2$ are isomorphic fields
and there exists a bicontinuous semi-linear bijection $S
:E_1\mapsto E_2$ such that, for every $T\in \mathcal{B}(E_1)$,
\begin{center}
$\Phi(T) = S TS\ep $.
\end{center}
\end{prop}
Proof. If $p$ is a continuous projection then, as projections are
defined by means of an equation in   the language of rings, $\Phi
(p)$ is also a continuous projection. The same argument shows
that, for two projections $p$ and $q$, we have $p\leq q$ if and
only if $\Phi(p) \leq \Phi(q)$. Moreover $\Phi(1_{E_1}-p)=1_{E_2}-
\Phi(p)$ and the restriction of $\Phi$ to the set of all
continuous projections is an orthoposet isomorphism.

 Fix a linear projection $p\in \mathcal{B}(E_1)$ of rank 1 (such
projection exists since the projection lattice of  the DAC-lattice
$L_{F_1}(E_1)$ is atomistic with atoms of the form $(X,Y)$, $X$ a
one dimensional subspace)   and non-zero elements $x_0\in Im\;p$,
$y_0\in Im\;\Phi(p)$. Remark that $\Phi(p)$ is also a continuous
projection of rank 1.

Let $x\in E_1$    and consider the linear mapping $U$ defined by
$U(x_0) = x$ and $U(t)= 0$ for $t\in Ker p$. The mapping $U$ is
continuous as a linear mapping with a finite-dimensional range and
a closed kernel (\cite{SCH},page 75).

Assume that $V\in\mathcal{B}(E_1)$ also satisfies $V(x_0)= x$. For
every $\lambda\in K_1$, $U(\lambda x_0) = V(\lambda x_0)$ and thus
$U\circ p = V\circ p$. We have $\Phi(U) \circ \Phi(p) = \Phi(V)
\circ \Phi(p)$ and so $\Phi(U) (y_0)=\Phi(V) (y_0)$. Thus, we can
define a  mapping $S: E_1\mapsto E_2$ by $S(x) = \Phi(U)(y_0)$.

Let $x,\;x\pr \in E_1$ and $U$, $U\pr$, $W$ be elements of
$\mathcal{B}(E_1 )$ such that $U(x_0)=x$, $U\pr(x_0) = x\pr$,
$W(x_0)=x+x\pr$. As $(U+U\pr)(x_0) =x+y$, we have $(U+U\pr)\circ
p=W\circ p$ and $S(x) +S(x\pr) = S(x+x\pr)$.

In a similar way, it can be proved that $S$ is a bijection and
$\Phi(T) = STS\ep$, for every $T\in\mathcal{B}(E_1)$.

The center of the rings $\mathcal{B}(E_i)$, $i=1,\, 2$ is
$\{k\,1_{E_i}\;|\; k\in K_i\}$ since a linear mapping which
commutes with all projections of rank 1 is a homothetic
transformation. Therefore,    for every $k\in K_1$, there exists $
k\pr \in K_2$ such that $\Phi(k 1_{E_1})= k\pr 1_{E_2}$. One can
check that the mapping $\sigma:K_1\mapsto K_2$ defined by
$\sigma(k)=k\pr$ is an isomorphism from the field $K_1$ onto the
field  $K_2$ and that $S(\lambda x) = \sigma( \lambda)S(x)$ for
every $\lambda\in K_1$ .

The last step is the proof of   continuity of $S$. Let $f\in F_1$
be a continuous non-zero linear form on $E_1$. Define a linear
mapping $T:E_1\mapsto E_1$ by $T(x) = f(x)x_0$. The mapping $x\in
E_1\mapsto (f(x),x_0)\in K\times E_1$ is continuous and, as $E_1$
is a topological vector space, $T\in \mathcal{B}(E_1)$. We have $S
TS\ep  = \Phi(T) \in\mathcal{B}(E_2)$ and
\begin{eqnarray*} x\in Ker\; STS\ep & \Leftrightarrow& TS(x)
=0\Leftrightarrow f(S  (x)) =0\\ &\Leftrightarrow&
 S (x) \in
Ker\;f\Leftrightarrow x\in S(Ker\; f)
\end{eqnarray*}
   Since $STS\ep= \Phi(T)$ is continuous, $Ker\; STS\ep $ is
closed and $S$, which carries hyperplanes to hyperplanes, carries
closed hyperplanes to closed hyperplanes. The mapping $S\ep$ is
continuous and, by symmetry, so does $S$.
\medskip

 \noindent{\bf Remark.} Assume that $E_1$ and $E_2$ are real locally convex spaces  and
that $F_1$ and $F_2$ are their topological duals for the weak
topology. The rings $ \mathcal{B}(E_1)$ and $\mathcal{B}(E_2)$ are
the rings of weakly continuous linear mappings defined on $E_1$
and $E_2$. Every ring isomorphism $ \Phi:\mathcal{B}(E_1) \mapsto
\mathcal{B}(E_2)$ is of the form $ \Phi (T)= STS\ep $ with $S:E_1
\mapsto E_2$ a weakly bicontinuous linear bijection. If $E_1$ and
$E_2$ are metrizable then $S$ is continuous (\cite{SCH}
chap.IV,7.4). For real Banach space this result is due to S.
Eidelheit (\cite{EI}).

Now assume that  $E_1$ and $E_2$ are infinite dimensional complex
normed spaces and $F_1 =E_1\pr$, $F_2=E_2\pr$. As $S$ carries
closed hyperplanes to closed hyperplanes, $S$ is linear or
conjugate linear (\cite{KA}, Lemma 2) and by \cite{FIL}, Lemma 3,
$S$ is bicontinuous. We have obtained a result of \cite{AR} (see
also \cite{FIL}, Theorem 2) : if $\Phi$ is an isomorphism of the
rings of continuous linear transformations on infinite dimensional
complex normed spaces $E_1$   and $E_2$ then there exists a
bicontinuous linear or conjugate linear bijection $S:E_1 \mapsto
E_2$ such that $\Phi(T) = STS\ep$.

\begin{prop} Let $(E,F)$ be a pair of dual
spaces.

\begin{enumerate}\item The restriction of an automorphism of the ring $ \mathcal{B}(E)
$ to the set of continuous linear projection is an  even
orthoposet automorphism and the restriction of an
anti-automorphism of the ring $ \mathcal{B}(E) $ to the set of
continuous linear projection is an  odd orthoposet automorphism.
 \item Conversely,    every even automorphism of the orthoposet of
continuous  linear projections defined on $E$ extends to an
automorphism of the ring $ \mathcal{B}(E)$.
\end{enumerate}
\end{prop}
Proof. 1) If $\phi$ is an automorphism or an anti-automorphism of
the ring $ \mathcal{B}(E)$ then its restriction  to the set of
continuous linear projections is an   orthoposet automorphism. The
nature of this automorphism will be given by the following lemma.
\begin{lem} Let $p$ and $q$ two linear projections defined on a
vector space $E$.
\begin{enumerate}
  \item $Im\, p =  Im\, q \Leftrightarrow pq = q$ and $qp=p$.
  \item $Ker\, p =Ker\, q \Leftrightarrow pq=p$ and $qp=q$.

\end{enumerate}
\end{lem}
Proof. If $pq=q$ and $x\in Im\, q$ then $q(x)=x$ and, since
$p(q(x))=q(x)$, we have $p(x) = x$ and $x\in Im\,p$. Conversely,
if $Im\,q \subset Im\, p$ and $x\in E$ then $x=x_1 +x_2$ with
$x_1\in Im\,q$ and $x_2\in Ker \, q$. We have $p(q(x))=
p(q(x_1))=p(x_1) = x_1 = q(x)$ and thus $p(q(x))=q(x)$. Finally,
$Im\,q \subset Im\, p \Leftrightarrow pq =q$ and all the other
proofs are similar. \medskip

We return to the proof of the proposition. If $\phi$
 is an automorphism of $ \mathcal{B}(E)$ then, for two projections
 $p$ and $q$, \begin{eqnarray*}Im\, p=Im\,q &\Leftrightarrow &pq=q\;\;
 \mathrm{et}\;\;qp= p\\&\Leftrightarrow &\phi(p)\phi(q)=\phi(q)\;\;
 \mathrm{et}\;\;\phi(q)\phi(p)= \phi(p)\\
& \Leftrightarrow& Im\,\phi(p)=Im\,\phi(q),
\end{eqnarray*}
 and the  restriction  of $\phi$ to the set of continuous
linear projections is an even  orthoposet automorphism. By a
similar proof, the  restriction    to the set of continuous linear
projections of an anti-automorphism is an odd  orthoposet
automorphism.

If $\Psi$ is an
 even orthoposet automorphism of the set of all continuous linear
 projection defined on $E$
 then there exists an automorphism $f$ of the lattice of all closed
 subspaces of $E$ such that $\Psi(Im\;p,Ker\;p) = (f(Im\;
p),f(Ker\; p))$. Let $S$ be  the bicontinuous semi-linear
bijection
  such that $S(X) =f(X)$ for every closed
subspace $X$ of $E$. We have $\Psi(p) = SpS\ep$ for every
projection $p\in \mathcal{B}(E)$ and $T \in \mathcal{B}(E) \mapsto
STS\ep$ is an automorphism of the ring $ \mathcal{B}(E)$ which
extends $\Psi$.
\bigskip

\noindent {\bf Question:} Do  odd automorphisms of the orthoposet
of continuous  linear projections defined on $E$ extend to
anti-automorphisms of the ring $ \mathcal{B}(E)$?

A problem in the description of odd automorphisms of lattices of
projections  is the lack of knowledge about anti-automorphisms of
lattices of closed subspaces in the infinite dimensional case.
Anti-automorphisms which are orthocomplementations are described,
as in the finite dimensional case, by means of symmetric bilinear
forms (\cite{VA}, Lemma 4.2.) but if an anti-automorphism $\Phi$
does not satisfy $M \cap \Phi(M)= \{0\}$ then, for a finite
dimensional subspace $F$, $M\mapsto \Phi(M) \cap F$ is not, in
general, an anti-automorphism of $[0,F]$ and it is not possible to
reduce the infinite dimensional case to the finite dimensional
one in the usual way.

 \end{document}